\documentclass[11pt]{article}
\usepackage{amssymb,amsmath,amsthm,pslatex,times}
\usepackage{subfigure,url,multirow}

\title{The metric dimension of small distance-regular and strongly regular graphs}
\author{Robert~F.~Bailey\footnote{Division of Science (Mathematics), Grenfell Campus, Memorial University of Newfoundland, University Drive, Corner Brook, NL A2H~6P9, Canada. E-mail: \texttt{rbailey@grenfell.mun.ca}}}

\newtheorem{thm}{Theorem}

\newtheorem{prop}[thm]{Proposition}

\theoremstyle{definition}

\renewcommand{\d}{\mathrm{d}}
\DeclareMathOperator{\Aut}{Aut}

\begin{document}

\maketitle

\begin{abstract}
A {\em resolving set} for a graph $\Gamma$ is a collection of vertices $S$, chosen so that for each vertex $v$, the list of distances from $v$ to the members of $S$ uniquely specifies $v$. The {\em metric dimension} of $\Gamma$ is the smallest size of a resolving set for $\Gamma$.  

A graph is {\em distance-regular} if, for any two vertices $u,v$ at each distance $i$, the number of neighbours of $v$ at each possible distance from $u$ (i.e.\ $i-1$, $i$ or $i+1$) depends only on the distance $i$, and not on the choice of vertices $u,v$.  The class of distance-regular graphs includes all distance-transitive graphs and all strongly regular graphs.

In this paper, we present the results of computer calculations which have found the metric dimension of all distance-regular graphs on up to 34 vertices, low-valency distance transitive graphs on up to 100 vertices, strongly regular graphs on up to 45 vertices, rank-$3$ strongly regular graphs on under 100 vertices, as well as certain other distance-regular graphs.
\\


\hrule
\end{abstract}

%\newpage

\section{Introduction}
%{\sf define met dim}\\
A {\em resolving set} for a graph $\Gamma=(V,E)$ is a set of vertices $S=\{v_1,\ldots,v_k\}$ such that for each vertex $w \in V$, the list of distances $(\d(w,v_1),\ldots,\d(w,v_k))$ uniquely determines $w$.  Equivalently, $S$ is a resolving set for $\Gamma$ if, for any pair of vertices $u,w \in V$, there exists $v_i \in S$ such that $\d(u,v_i)\neq \d(w,v_i)$; we say that $v_i$ {\em resolves} $u$ and $w$.  The {\em metric dimension} of $\Gamma$ is the smallest size of a resolving set for $\Gamma$.  This concept was introduced to the graph theory literature in the 1970s by Harary and Melter~\cite{Harary76} and, independently, Slater~\cite{Slater75}; however, in the context of arbitrary metric spaces, the concept dates back at least as far as the 1950s (see Blumenthal~\cite{Blumenthal53}, for instance).  In more recent years, there has been a considerable number of papers written about the metric dimension of graphs.  For further details, the reader is referred to the survey paper \cite{bsmd}.

%{\sf define DRG, SRG, AS, PPoly}\\
When studying metric dimension, distance-regular graphs are a natural class of graphs to consider.  A graph $\Gamma$ is {\em distance-regular} if, for all $i$ with $0\leq i\leq \mathrm{diam}(\Gamma)$ and any vertices $u,v$ with $\d(u,v)=i$, the number of neighbours of $v$ at distances $i-1$, $i$ and $i+1$ from $u$ depend only on the distance $i$, and not on the choices of $u$ and $v$.  For more information about distance-regular graphs, see the book of Brouwer, Cohen and Neumaier~\cite{BCN} and the forthcoming survey paper by van Dam, Koolen and Tanaka~\cite{vanDam}.  Note that the class of distance-regular graphs contains the distance-transitive graphs (i.e.\ those graphs $\Gamma$ with the property that for any vertices $u,v,u',v'$ such that $\d(u,v)=\d(u',v')$, there exists an automorphism $g$ such that $u^g=u'$ and $v^g=v'$) and the connected strongly regular graphs (which are the distance-regular graphs of diameter~$2$).  

For any graph $\Gamma$ with diameter~$d$, consider the partition of $V\times V$ into $d+1$ parts, given by the pairs of vertices at each possible distance in $\Gamma$.  If $\Gamma$ is distance-regular, this partition is an example of an {\em association scheme}.  These are much more general objects, and ones which are inconsistently named in the literature; see \cite[Section 3.3]{bsmd} for more details.  An association scheme is said to be {\em $P$-polynomial} if it arises from a distance-regular graph; however, more than one graph may give rise to the same $P$-polynomial association scheme (see \cite{BCN,vanDam}).  It is not difficult to see that two graphs arising in this way must have the same metric dimension (see \cite[Section 3.5]{bsmd}).

Since the publication of the survey paper by Cameron and the present author~\cite{bsmd}, a number of papers have been written on determining, or bounding, the metric dimension of various families of distance-regular graphs: see \cite{imprimitive,jk,grassmann,Beardon,FengWang,GuoWangLi,GuoWangLi2,Heger}, for instance.  The purpose of this paper is to give the results of a number of computer calculations, using the {\sf GAP} computer algebra system~\cite{gap}, which have obtained the metric dimension for all ``small'' distance-regular graphs (i.e.\ on up to 34 vertices), for distance-regular graphs of valency 3 and 4 (on up to 189 vertices), low-valency distance-transitive graphs (up to valency 13, and up to 100 vertices), and certain other distance-regular graphs.  For strongly regular graphs, this has included an independent verification of earlier computations by Kratica {\em et al.}~\cite{Kratica08} (which used an entirely different approach via linear programming).

\section{Known results}
%{\sf state previous results}\\
%{\sf complete, complete bipartite, cycles, Hamming H(2,n), Johnson J(n,2),.....}
In this section, we summarize the relevant known values of the metric dimension of various families of distance-regular graphs.  It is a straightforward exercise to verify that the complete graph $K_n$ has metric dimension $n-1$, that the complete bipartite graph $K_{m,n}$ has metric dimension $m+n-2$, and that a cycle $C_n$ with $n\geq 3$ vertices has metric dimension~$2$.  The following result is also straightforward, yet the author is not aware of it appearing anywhere in the literature.

\begin{prop} \label{prop:multipartite}
%Consider a complete multipartite graph $K_{m^r}$ with $r$ parts of size $m$, for $r>1$.  Then the metric dimension of $K_{m^r}$ is $r(m-1)$.
Consider a complete multipartite graph $\Gamma=K_{m_1,\ldots,m_r}$ with $r$ parts of sizes $m_1,\ldots,m_r$, for $r>1$.  Then the metric dimension of $\Gamma$ is $\sum_{i=1}^r (m_i -1)$.
\end{prop}

In particular, in the special case where a complete multipartite graph has $r$ parts of size $m$ (and thus is strongly regular), %Proposition~\ref{prop:multipartite} 
this shows that the metric dimension is $r(m-1)$.

\proof Suppose the vertex set of $\Gamma$ is $V=V_1 \cup \cdots \cup V_r$, where the $V_i$ are disjoint sets of sizes $m_1,\ldots,m_r$; every possible edge exists from $V_i$ to $V_j$ (for $i\neq j$), and no edges exist inside any $V_i$.  Let $T$ 
%$=\{w_1,\ldots,w_r\}$
be a transversal of $V_1,\ldots,V_r$, and let $S=V\setminus T$.  It is straightforward to verify that $S$ is a resolving set for $\Gamma$ of size $\sum_{i=1}^r (m_i -1)$.  Furthermore, no smaller resolving set may exist: suppose for a contradiction that $R$ is a subset of $V$ with size smaller than the above.  By the pigeonhole principle there exists an index $i$ for which $V_i$ contains two vertices $u,v$ not in $R$, and no vertex in $R$ will resolve this pair of vertices. \endproof

The following results are somewhat less trivial.  Recall that the {\em Johnson graph} $J(n,k)$ is the graph whose vertex set consists of all $k$-subsets of an $n$-set, and two $k$-subsets are adjacent if and only if they intersect in a $(k-1)$-subset.  The {\em Kneser graph} $K(n,k)$ has the same vertex set as $J(n,k)$, but adjacency is defined by two $k$-sets being disjoint.  The Johnson graph is always distance-regular, whereas the Kneser graph only is in two special cases, namely $K(n,2)$ (which is the complement of $J(n,2)$) and $K(2k+1,k)$ (known as the {\em Odd graph}, and usually denoted $O_{k+1}$). The following result was obtained by Cameron and the present author in 2011.

\begin{thm}[%  %%%DON'T REMOVE THE '%' AND LINEBREAKS HERE--IT CONFUSES THE LaTeX COMPILER!!!
{Bailey and Cameron \cite[Corollary 3.33]{bsmd}}%
] 
\label{thm:johnson}
For $n \geq 6$, metric dimension of the \linebreak Johnson graph $J(n,2)$ and Kneser graph $K(n,2)$ is $\frac{2}{3}(n-i)+i$, where $n\equiv i \pmod 3$,\linebreak $i\in \{0,1,2\}$.
\end{thm}

It is easy to determine the metric dimension of $J(3,2)$, $J(4,2)$ and $J(5,2)$ by hand: these values are $2$, $3$ and $3$, respectively.  Further results about resolving sets for Johnson and Kneser graphs may be found in~\cite{jk}.

We also recall that the {\em Hamming graph} $H(d,q)$ has as its vertex set the collection of all $d$-tuples over an alphabet of size $q$, and two $d$-tuples are adjacent if and only if they differ in exactly one position; these graphs are distance-transitive.  Two important examples are the hypercube $H(d,2)$ and the square lattice graph $H(2,q)$.  The following result was obtained by C\'aceres {\em et al}.\ in 2007.

\begin{thm}[%
{C\'aceres {\em et al.}\ \cite[Theorem 6.1]{cartesian}}%
] 
\label{thm:lattice}
For all $q \geq 1$, the metric dimension of the square lattice graph $H(2,q)$ is $\lfloor \frac{2}{3}(2q-1) \rfloor$.
\end{thm}

This is the only infinite family of Hamming graphs for which the metric dimension is known precisely.  Further details about the metric dimension of Hamming graphs can be found in \cite[Section 3.6]{bsmd}; for the hypercubes, see also Beardon~\cite{Beardon}.  Some precise values were computed by Kratica {\em et al}.\ \cite{Kratica09}, using genetic algorithms.

Finally, we mention a recent result of the present author, which helps to eliminate the need for some additional computations.  The {\em bipartite double} of a graph $\Gamma=(V,E)$ is a bipartite graph $D(\Gamma)$, whose vertex set consists of two disjoint copies of $V$, labelled $V^+$ and $V^-$, with $v^+$ adjacent to $w^-$ in $D(\Gamma)$ if $v$ and $w$ are adjacent in $\Gamma$.

\begin{thm}[%
{Bailey \cite{imprimitive}}%
]
\label{thm:double}
Suppose $\Gamma$ is a distance-regular graph of diameter $d$, and whose shortest odd cycle has length $2d+1$.  Then $\Gamma$ and its bipartite double $D(\Gamma)$ have the same metric dimension.
\end{thm}

In particular, the graph $K_{n,n}-I$ (obtained by deleting a perfect matching from $K_{n,n}$) is the bipartite double of the complete graph $K_n$, which satisfies the conditions of Theorem~\ref{thm:double}, and thus $K_{n,n}-I$ has metric dimension $n-1$.

%\newpage

\section{The method}
Suppose $\Gamma=(V,E)$ is a graph with $n$ vertices, labelled $v_1,\ldots,v_n$.  The {\em distance matrix} of $\Gamma$ is the $n\times n$ matrix $A$, whose $(i,j)$ entry is the distance in $\Gamma$ from $v_i$ to $v_j$.  %In the case where $\Gamma$ is distance-regular, this is equivalent to the class matrix of the corresponding association scheme.  
Note that if there are two distance-regular graphs arising from the same $P$-polynomial association scheme, their distance matrices are equivalent, up to a relabelling of the distance classes.  For instance, the metric dimension of a primitive strongly regular graph and that of its complement will be equal.

Suppose $A$ is the distance matrix of $\Gamma$, and let $S$ be a subset of $V$.  Denote by $[A]_S$ the submatrix formed by taking the columns of $A$ indexed by elements of $S$.  It is clear that $S$ is a resolving set for $\Gamma$ if and only if the rows of $[A]_S$ are distinct.  Furthermore, if one has obtained the distance matrix of a relatively small graph, it is near-instantaneous for a computer to verify if a given submatrix has this property.

The results given below were obtained using the computer algebra system ${\sf GAP}$ \cite{gap} and various packages developed for it, in particular the {\sf GRAPE} package of Soicher \cite{grape}, the functions for association schemes of Hanaki \cite{Hanaki-gap}, and also some functions of Cameron \cite{cameron-gap}.  However, the most useful tool for these computations is the {\sf SetOrbit} package of Pech and Reichard \cite{setorbit}.  Given a set $V$ and a group $G$ acting on it, this provides an efficient method for enumerating a canonical representative of each orbit of $G$ acting on the subsets of $V$ of a given size.  It is clear that $S$ is a resolving set for a graph $\Gamma$ if and only if its image $S^g = \{x^g \, : \, x\in S\}$ is a resolving set, for any $g\in \Aut(\Gamma)$.  Consequently, when searching for a resolving set of a particular size, it suffices to test just one representative of each orbit on subsets of that size.  Therefore, the methods of Pech and Reichard (which are explained in detail in~\cite{PechReichard}) are precisely what is needed to dramatically reduce the search space.

\subsection{Data sources}
Hanaki and Miyamoto's library of small association schemes \cite{Hanaki-data} contains all distance-regular graphs on up to 34 vertices: one merely has to filter out the $P$-polynomial examples from their lists.  Association schemes with primitive automorphism groups may be constructed by using {\sf GAP}'s internal libraries of primitive groups.  Graphs (such as incidence graphs and point graphs) obtained from generalized polygons may be constructed using the {\sf FinInG} package of Bamberg {\em et al.}\ \cite{fining}, while for graphs associated with block designs (including projective and affine geometries) the {\sf DESIGN} package of Soicher \cite{design} is used. Other useful data sources include Spence's catalogue of strongly regular graphs \cite{spence}, Royle's catalogue of symmetric (i.e.\ arc-transitive) cubic graphs \cite{Royle}, the online \textsc{Atlas} of Finite Group Representations \cite{wwwATLAS}, and Sloane's libraries of Hadamard matrices \cite{sloane}.

%{\sf uses GAP, GRAPE, FinInG, DESIGN, Association Schemes, SetOrbit}, 

%{\sf sources: Hanaki, Spence, wwwATLAS, Sloane, GAP libs}

%\newpage

\section{Results}
\renewcommand{\arraystretch}{1.1}
\subsection{Distance-regular graphs on up to 34 vertices} \label{subsection:smalldrg}
As mentioned above, all distance-regular graphs on up to 34 may all be obtained from the catalogue of association schemes of Hanaki and Miyamoto~\cite{Hanaki-data}.  For 31 and 32 vertices, a complete classification of association schemes is not available; however, as there are no feasible parameter sets for strongly regular graphs with these numbers of vertices, and (other than $K_{31}$ and $C_{31}$) no distance-transitive graphs with 31 vertices, one can see that the data available contains all distance-regular graphs with 31 and 32 vertices.  These graphs are all small enough that no sophisticated computations were required to determine the metric dimension; the results are given in Tables \ref{table:smalldrg1} and \ref{table:smalldrg2}.  Cycles, complete graphs, complete bipartite graphs and complete multipartite graphs are omitted from these tables.  Those graphs which are not distance-transitive are indicated $\dagger$; cases where there is a unique distance-transitive example are indicated $\ddag$.  The abbreviation IG is used to denote an incidence graph.

\begin{table}[p]
  \centering
  \begin{tabular}{|c|l|c|c|c|}
  \hline
  No.\ of vertices & Graph					     & Valency & Diameter & Met.\ dim. \\
  \hline
  6    & Octahedron $J(4,2)$    				 & 4       & 2        & 3 \\
  \hline
  8    & Cube $Q_3 \cong K_{4,4}-I$			 & 3       & 3        & 3 \\
  \hline
  9    & Paley graph $P_9 \cong H(2,3)$	 & 4       & 2        & 3 \\
  \hline
  \multirow{3}{*}{10}
       & Petersen graph $O_3 = K(5,2)$	 & 3       & 2        & 3 \\
       & $J(5,2)$     									 & 6       & 2        & 3 \\
       & $K_{5,5}-I$                     & 4       & 3        & 4 \\
  \hline
  \multirow{2}{*}{12}
       & Icosahedron                     & 5       & 3        & 3 \\
       & $K_{6,6}-I$                     & 5       & 3        & 5 \\
  \hline
  13   & Paley graph $P_{13}$	        	 & 6       & 2        & 4 \\
  \hline
  \multirow{3}{*}{14}
       & Heawood (IG of $\mathrm{PG}(2,2)$) 
                                         & 3       & 3        & 5 \\
       & Non-IG of $\mathrm{PG}(2,2)$    & 4       & 3        & 5 \\
       & $K_{7,7}-I$                     & 6       & 3        & 6 \\
  \hline
  \multirow{3}{*}{15}
       & Line graph of Petersen graph    & 4       & 3        & 4 \\
       & $K(6,2)$                        & 6       & 2        & 4 \\
       & $J(6,2)$                        & 8       & 2        & 4 \\
  \hline
  \multirow{8}{*}{16}
       & $4$-cube $Q_4$                  & 4       & 4        & 4 \\
       & $H(2,4)$                        & 6       & 2        & 4 \\
       & Complement of $H(2,4)$          & 9       & 2        & 4 \\
       & Shrikhande graph$\dagger$       & 6       & 2        & 4 \\
       & Complement of Shrikhande graph$\dagger$
                                         & 9       & 2        & 4 \\
       & Clebsch graph                   & 5       & 2        & 4 \\
       & (Complement of) Clebsch graph   & 10      & 2        & 4 \\
       & $K_{8,8}-I$                     & 7       & 3        & 7 \\
  \hline
  17   & Paley graph $P_{17}$	        	 & 8       & 2        & 4 \\
  \hline
  \multirow{2}{*}{18}
       & Pappus graph                    & 3       & 4        & 4 \\
       & $K_{9,9}-I$                     & 8       & 3        & 8 \\
  \hline
  \multirow{4}{*}{20}
       & Dodecahedron                    & 3       & 5        & 3 \\
       & Desargues graph $D(O_3)$        & 3       & 5        & 3 \\
       & $J(6,3)$                        & 9       & 3        & 4 \\       
       & $K_{10,10}-I$                   & 9       & 3        & 9 \\
  \hline
  \multirow{3}{*}{21}
       & Line graph of Heawood graph     & 4       & 3        & 4 \\
       & $J(7,2)$                        & 10      & 2        & 4 \\
       & $K(7,2)$                        & 10      & 2        & 4 \\
  \hline
  \multirow{3}{*}{22}
       & IG of biplane                   & 5       & 3        & 6 \\
       & Non-IG of biplane               & 6       & 3        & 6 \\
       & $K_{11,11}-I$                   & 10      & 3        & 10 \\
  \hline
  \multirow{2}{*}{24}
       & Symplectic cover \cite[p.\ 386]{BCN}$\dagger$
                                         & 7       & 3        & 5 \\
       & $K_{12,12}-I$                   & 11      & 3        & 11 \\
  \hline
  \end{tabular}
\caption{Metric dimension of distance-regular graphs on up to 24 vertices \label{table:smalldrg1}}
\end{table}

%\newpage
  
\begin{table}[p]
  \centering
  \begin{tabular}{|c|l|c|c|c|}
  \hline
  No.\ of vertices & Graph					     & Valency & Diameter & Met.\ dim. \\
  \hline
  \multirow{4}{*}{25}
       & $H(2,5)$                        & 8       & 2        & 6 \\
       & Paley graph $P_{25}$            & 12      & 2        & 5 \\
       & Other $\mathrm{srg}(25,12,5,6)$ {\em ($14$ graphs, $7$ pairs)}$\dagger$
                                         & 12      & 2        & 5 \\
       & Complement of $H(2,5)$          & 16      & 2        & 6 \\
  \hline
  \multirow{5}{*}{26}
       & $\mathrm{srg}(26,10,3,4)$ {\em ($10$ graphs)}$\dagger$
                                         & 10      & 2        & 5 \\
       & IG of $\mathrm{PG}(2,3)$        & 4       & 3        & 8 \\
       & Non-IG of $\mathrm{PG}(2,3)$    & 9       & 3        & 8 \\
       & $K_{13,13}-I$                   & 12      & 3        & 12 \\
       & Complements of $\mathrm{srg}(26,10,3,4)$ {\em ($10$ graphs)}$\dagger$
                                         & 15      & 2        & 5 \\
  \hline
  \multirow{4}{*}{27}
       & $H(3,3)$                        & 6       & 3        & 4 \\
       & $GQ(2,4)$ minus spread {\em ($2$ graphs)}
                                         & 8       & 3        & 5 \\
       & Complement of Schl\"afi graph   & 10      & 2        & 5 \\
       & Schl\"afi graph                 & 10      & 2        & 5 \\
  \hline
  \multirow{7}{*}{28}
			 & Coxeter graph									 & 3       & 4        & 4 \\
			 & $J(8,2)$                        & 12      & 2        & 6 \\
			 & Chang graphs {\em ($3$ graphs)}$\dagger$
			                                   & 12      & 2        & 6 \\
			 & Taylor graph from $P_{13}$      & 13      & 3        & 5 \\
             & $K_{14,14}-I$                   & 13      & 3        & 13 \\
			 & $K(8,2)$                        & 15      & 2        & 6 \\
			 & Complements of Chang graphs {\em ($3$ graphs)}$\dagger$
			                                   & 15      & 2        & 6 \\
  \hline
  \multirow{2}{*}{29}
       & Paley graph $P_{29}$            & 14      & 2        & 6 \\
       & Other $\mathrm{srg}(29,14,6,7)$ {\em ($40$ graphs, $20$ pairs)}$\dagger$
                                         & 14      & 2        & 5 \\
  \hline
  \multirow{6}{*}{30}
       & Tutte's 8-cage                  & 3       & 4        & 6 \\
       & IG of $\mathrm{PG}(3,2)$        & 7       & 3        & 8 \\
       & Non-IG of $\mathrm{PG}(3,2)$    & 8       & 3        & 8 \\
       & IGs of Hadamard designs {\em ($3$ graphs)}$\ddag$
                                         & 7       & 3        & 8 \\
       & Non-IGs of Hadamard designs {\em ($3$ graphs)}$\ddag$
                                         & 8       & 3        & 8 \\
       & $K_{15,15}-I$                   & 14      & 3        & 14 \\
  \hline
  \multirow{9}{*}{32}
       & IG of truncated $\mathrm{AG}(2,4)$
                                         & 4       & 4        & 6 \\
       & 5-cube $Q_5$                    & 5       & 5        & 4 \\
       & Armanios--Wells graph           & 5       & 4        & 5 \\
       & IGs of biplanes {\em ($3$ graphs)}$\dagger$
                                         & 6       & 3        & 8 \\
       & Hadamard graph                  & 8       & 4        & 7 \\
       & Non-IGs of biplanes {\em ($3$ graphs)}$\dagger$
                                         & 10      & 3        & 8 \\
       & Taylor graph from $J(6,2)$      & 15      & 3        & 5 \\
       & Taylor graph from $K(6,2)$      & 15      & 3        & 5 \\
       & $K_{16,16}-I$                   & 15      & 3        & 15 \\
  \hline
  34   & $K_{17,17}-I$                   & 16      & 3        & 16 \\
  \hline
  \end{tabular}
\caption{Metric dimension of distance-regular graphs on 25 to 34 vertices \label{table:smalldrg2}}
\end{table}

%\newpage

\subsection{Distance-regular graphs of valency 3 and 4} \label{subsection:cubic}
The distance-transitive graphs of valency 3 were determined by Biggs and Smith in 1971 \cite{BiggsSmith71} (see also Gardiner \cite{Gardiner75}); this classification was extended to all distance-regular graphs of valency 3 by Biggs, Boshier and Shawe-Taylor in 1986 \cite{Biggs86}.  In addition to $K_4$ and $K_{3,3}$, there are eleven such distance-regular graphs, of which all but one are distance-transitive, with the exception being Tutte's 12-cage.  The metric dimension of each of these graphs is given in Table~\ref{table:valency3}.

The distance-transitive graphs of valency~4 were determined in 1974 by Smith \cite{Smith73,Smith74a,Smith74b} (see also Gardiner \cite{Gardiner85}); this classification was extended to distance-regular graphs by Brouwer and Koolen in 1999 \cite{BrouwerKoolen}.  The classification is complete, except for relying on a classification of generalized hexagons $GH(3,3)$, which give rise to 4-regular distance-regular graphs on 728 vertices.  (In any case, these would be beyond the scope of the computations in this paper.)  The metric dimension of each of these graphs is given in Table~\ref{table:valency4}.

%Most of these examples can be obtained using 
Of the graphs with more than 34 vertices, the Foster graph was obtained from Royle's \linebreak catalogue \cite{Royle}, the Biggs--Smith graph from its automorphism group $\mathrm{PSL}(2,17)$, and Tutte's \mbox{12-cage} as the incidence graph of the generalized hexagon $GH(2,2)$ (constructed using \mbox{{\sf FinInG}).}  
The Odd graph $O_4$ lies inside the Johnson scheme $J(7,3)$, which may be constructed using {\sf GRAPE}, as can line graphs and bipartite doubles.  The incidence graph of the $GQ(3,3)$ may be constructed using {\sf FinInG}.

\begin{table}[hbtp]
  \centering
  \begin{tabular}{|l|c|c|c|}
  \hline
  Graph & No.\ of vertices & Diameter & Metric dimension \\
  \hline
%  $K_4$                       & 4 & 1 & 3 \\
%  $K_{3,3}$                   & 6 & 2 & 4 \\
  Cube $Q_3 \cong K_{4,4}-I$  & 8 & 3 & 3 \\
  Petersen graph $O_3$        &10 & 2 & 3 \\
  Heawood graph               &14 & 3 & 5 \\
  Pappus graph                &18 & 4 & 4 \\
  Dodecahedron                &20 & 5 & 3 \\
  Desargues graph $D(O_3)$    &20 & 5 & 3 \\
  Coxeter graph               &28 & 4 & 4 \\
  Tutte's 8-cage              &30 & 4 & 6 \\
  Foster graph                &90 & 8 & 5 \\
  Biggs--Smith graph         &102 & 7 & 4 \\
  Tutte's 12-cage$\dagger$   &126 & 6 & 8 \\ \hline
  \end{tabular}
\caption{Metric dimension of distance-regular graphs of valency 3 \label{table:valency3}}
\end{table}

\begin{table}[hbtp]
  \centering
  \begin{tabular}{|l|c|c|c|}
  \hline
  Graph & No.\ of vertices & Diameter & Metric dimension \\
  \hline
%  $K_5$                       & 5 & 1 &4 \\
Octahedron $J(4,2)$ & 6        & 2       & 3 \\%       K_{2^3}
%K_{4,4}           8         2        6
Paley graph $P_9 \cong H(2,3)$
                   & 9        & 2       & 3 \\
$K_{5,5}-I$        & 10       & 3       & 5 \\
Distance-3 graph of Heawood graph
                   & 14       & 3       & 5 \\%       Non-incidence graph PG(2,2)
Line graph of Petersen graph
                   & 15       & 3       & 4 \\%       L(Petersen)
4-cube $Q_4$       & 16       & 4       & 4 \\
Line graph of Heawood graph
                   & 21       & 3       & 4 \\
Incidence graph of $\mathrm{PG}(2,3)$
                   & 26       & 3       & 8 \\%       mu=4q-4 (Bill)
Incidence graph of truncated $\mathrm{AG}(2,4)$
                   & 32       & 4       & 6 \\%       4-fold covering of K_{4,4}
Odd graph $O_4$    & 35       & 3       & 5 \\%       \leq 6 in BCGGMMP
Line graph of Tutte's 8-cage
                   & 45       & 4       & 4 \\
Doubled Odd graph $D(O_4)$
                   & 70       & 7       & 5 \\
Incidence graph of $GQ(3,3)${$\dagger$}
                   & 80       & 4       & 10 \\
Line graph of Tutte's 12-cage$\dagger$
                   & 189      & 6       & 6 \\ 
Incidence graph of $GH(3,3)$
                   & 728      & 6       & unknown \\ \hline
  \end{tabular}
\caption{Metric dimension of distance-regular graphs of valency 4 \label{table:valency4}}
\end{table}

\subsection{Distance-transitive graphs of valencies 5 to 13} \label{subsection:lowvalency}

The distance-transitive graphs of valencies 5, 6 and 7 were determined in 1986 by Faradjev, Ivanov and Ivanov \cite{Faradjev86}, and independently (for valencies 5 and 6) by Gardiner and Praeger \cite{GardinerPraeger86,GardinerPraeger87}.  For valencies 8 to 13, the classification was obtained in 1988 by Ivanov and Ivanov \cite{Ivanov88}.  

Tables \ref{table:valency5} to \ref{table:valency13} contain the metric dimension of all such graphs on up to 100 vertices (apart one exception on 98 vertices), as well as larger graphs when the computations succeeded.  No table is provided for valency~$11$, as the only example under 100 vertices (other than $K_{12}$ or $K_{11,11}$) is $K_{12,12}-I$, which has metric dimension~$11$.  An asterisk indicates that the metric dimension was not computed directly, but rather that Theorem~\ref{thm:double} was applied to the result of an earlier computation.

Of the examples with more than 34 vertices, many of the graphs are incidence graphs of designs or geometries, so can be constructed using the {\sf DESIGN} or {\sf FinInG} packages; the resolvable transversal designs $\mathrm{RT}[8,2;4]$ and $\mathrm{RT}[9,3;3]$ are given in the paper of Hanani~\cite{Hanani74}.  Otherwise, graphs were constructed in {\sf GAP} from their automorphism groups, either from the internal libraries of primitive groups, or %(such as with the Conway-Smith graph) 
using permutation representations in the \textsc{Atlas}.

\begin{table}[hbtp]
  \centering
  \begin{tabular}{|l|c|c|c|}
  \hline
  Graph & No.\ of vertices & Diameter & Metric dimension \\
  \hline
$K_{6,6}-I$       & 12     & 3        & 5\\
Icosahedron       & 12     & 3        & 3\\
Clebsch graph     & 16     & 2        & 4\\%        Folded 5-cube
Incidence graph of biplane
                  & 22     & 3        & 6\\
5-cube $Q_5$      & 32     & 5        & 4\\
Armanios--Wells graph
                  & 32     & 4        & 5\\
Sylvester graph from $\Aut(S_6)$
                  & 36     & 3        & 5\\
Incidence graph of $\mathrm{PG}(2,4)$
                  & 42     & 3        & 10\\%       IG PG(2,4): beats Bill's construction
Incidence graph of truncated $\mathrm{AG}(2,5)$
                  & 50     & 4        & 9\\
Odd graph $O_5$   & 126    & 4        & 6\\%       BCGGMMP \leq 8
Incidence graph of $GQ(4,4)$
                  & 170    & 4        & unknown\\
Doubled Odd graph $D(O_5)$
                  & 252    & 9        & $6^\ast$\\
  \hline
  \end{tabular}
\caption{Metric dimension of distance-transitive graphs of valency 5 \label{table:valency5}}
\end{table}

\begin{table}[hbtp]
  \centering
  \begin{tabular}{|l|c|c|c|}
  \hline
  Graph & No.\ of vertices & Diameter & Metric dimension \\
  \hline
%K_7               7         1        6
%K_{6,6}          12         2        10
$K_{7,7}-I$       & 14     & 3        & 6\\
$J(5,2)$          & 10     & 2        & 3\\%         co-Petersen
Paley graph $P_{13}$
                  & 13     & 2        & 4\\
$K(6,2)$          & 15     & 2        & 4\\
$H(2,4)$          & 16     & 2        & 4\\%         Spanish Cartesian paper
Non-incidence graph of biplane
                  & 22     & 3        & 6\\
$H(3,3)$          & 27     & 3        & 4\\
Incidence graph of biplane
                  & 32     & 3        & 8\\%         folded 6-cube
Hexacode graph    & 36     & 4        & 7\\
$2^{\textnormal{nd}}$ subconstituent of Hoffman--Singleton %graph
                 & 42      & 3        & 7\\
Halved Foster graph
                 & 45      & 4        & 6\\
Flag graph of $\mathrm{PG}(2,3)$
                 & 52      & 3        & 6\\
Perkel graph     & 57      & 3        & 6\\
Incidence graph of $\mathrm{PG}(2,5)$
                 & 62      & 3        & 15\\%       took 24.5 days.....
Point graph of $GH(2,2)$
                 & 63      & 3        & 6\\
Point graph of dual $GH(2,2)$
                 & 63      & 3        & 6\\
6-cube $Q_6$     & 64      & 6        & 5\\
  \hline
  \end{tabular}
\caption{Metric dimension of distance-transitive graphs of valency 6 \label{table:valency6}}
\end{table}

\begin{table}[hbtp]
  \centering
  \begin{tabular}{|l|c|c|c|}
  \hline
  Graph & No.\ of vertices & Diameter & Metric dimension \\
  \hline
$K_{8,8}-I$       & 16     & 3        & 7\\
Incidence graph of $\mathrm{PG}(3,2)$
                  & 30     & 3        & 8\\
Hoffman--Singleton graph
                  & 50     & 2        & 11\\
Folded 7-cube     & 64     & 3        & 6\\
Incidence graph of truncated $\mathrm{AG}(2,7)$
                  & 98     & 4        & unknown\\
Doubled Hoffman--Singleton graph
                  & 100    & 5        & $11^\ast$\\%  bipartite double of max odd girth     
7-cube $Q_7$      & 128    & 7        & 6\\
  \hline
  \end{tabular}
\caption{Metric dimension of distance-transitive graphs of valency 7 \label{table:valency7}}
\end{table}

\begin{table}[hbtp]
  \centering
  \begin{tabular}{|l|c|c|c|}
  \hline
  Graph & No.\ of vertices & Diameter & Metric dimension \\
  \hline
$J(6,2)$          & 15	   & 2        & 4\\
Paley graph $P_{17}$
                  & 17	   & 2        & 4\\
$K_{9,9}-I$       & 18     & 3        & 8\\
$H(2,5)$          & 25     & 2        & 6\\
Point graph of $GQ(2,4)$ minus spread
                  & 27     & 3        & 5\\        %2 graphs?
Non-incidence graph of $\mathrm{PG}(3,2)$
                  & 30     & 3        & 8\\
Hadamard graph	  & 32     & 4        & 7\\
Incidence graph of $\mathrm{RT}[8,2;4]$
                  & 64     & 4        & 10\\
$H(4,3)$		  & 81     & 4        & 5\\
Flag graph of $\mathrm{PG}(2,4)$
                  & 105    & 3        & 7\\
Folded $8$-cube   & 128    & 4        & 11\\
8-cube $Q_8$      & 256    & 8        & 6\\
  \hline
  \end{tabular}
\caption{Metric dimension of distance-transitive graphs of valency 8 \label{table:valency8}}
\end{table}

\begin{table}[hbtp]
  \centering
  \begin{tabular}{|l|c|c|c|}
  \hline
  Graph & No.\ of vertices & Diameter & Metric dimension \\
  \hline
Complement of $H(2,4)$
                    & 16     & 2        & 4\\
$K_{10,10}-I$       & 20     & 3        & 9\\
$J(6,3)$            & 20     & 3        & 4\\
Non-incidence graph of $\mathrm{PG}(2,3)$
                    & 26     & 3        & 8\\
Incidence graph of $\mathrm{RT}[9,3;3]$
                    & 54     & 4        & 10\\
$H(3,4)$            & 64     & 3        & 6\\
Unitals in $\mathrm{PG}(2,4)$ (from $\mathrm{PSL}(3,4)$)
                    & 280    & 4        & 5\\
9-cube $Q_9$        & 512    & 9        & 7\\
  \hline
  \end{tabular}
\caption{Metric dimension of distance-transitive graphs of valency 9 \label{table:valency9}}
\end{table}

\begin{table}[hbtp]
  \centering
  \begin{tabular}{|l|c|c|c|}
  \hline
  Graph & No.\ of vertices & Diameter & Metric dimension \\
  \hline
Clebsch graph       & 16     & 2        & 4\\
$J(7,2)$, $K(7,2)$  & 21     & 2        & 4\\
$K_{11,11}-I$       & 22     & 3        & 10\\
Complement of Schl\"afi graph
                    & 27     & 2        & 5\\
Non-incidence graph of biplane
                    & 32     & 3        & 8\\       %also 2 non-DT DRGs
$H(2,6)$            & 36     & 2        & 7\\
Gewirtz graph       & 56     & 2        & 9\\
Conway-Smith graph from $3.S_7$
                    & 63     & 4        & 6\\
Hall graph from $\mathrm{P\Sigma L}(2,25)$
                    & 65     & 3        & 6\\

Doubled Gewirtz graph
                    & 112    & 5        & $9^\ast$\\
$H(5,3)$            & 243    & 5        & 5\\
Hall--Janko near octagon from $\mathrm{J}_2.2$
                    & 315    & 4        & 8\\
  \hline
  \end{tabular}
\caption{Metric dimension of distance-transitive graphs of valency 10 \label{table:valency10}}
\end{table}

\begin{table}[hbtp]
  \centering
  \begin{tabular}{|l|c|c|c|}
  \hline
  Graph & No.\ of vertices & Diameter & Metric dimension \\
  \hline
Paley graph $P_{25}$ & 25    & 2        & 5\\
$K_{13,13}-I$       & 26     & 3        & 12\\
$J(8,2)$         	& 28     & 2        & 6\\
$J(7,3)$            & 35     & 3        & 5\\
Point graph of $GQ(3,3)$
                    & 40     & 2        & 7\\%         srg(40,12,2,4)
Point graph of dual $GQ(3,3)$
                    & 40     & 2        & 8\\%         srg(40,12,2,4)
Point graph of $GQ(4,2)$
                    & 45     & 2        & 8\\%         srg(45,12,3,3)
Hadamard graph      & 48     & 4        & 8\\
$H(2,7)$            & 49     & 2        & 8\\
Doro graph from $\mathrm{P\Sigma L}(2,16)$
                    & 68     & 3        & 6\\
$H(3,5)$            & 125    & 3        & 7\\
$H(4,4)$            & 256    & 4        & 7\\
  \hline
  \end{tabular}
\caption{Metric dimension of distance-transitive graphs of valency 12 \label{table:valency12}}
\end{table}

\begin{table}[hbtp]
  \centering
  \begin{tabular}{|l|c|c|c|}
  \hline
  Graph & No.\ of vertices & Diameter & Metric dimension \\
  \hline
$K_{14,14}-I$       & 28     & 3        & 13\\
Taylor graph from $P_{13}$
                    & 28     & 3        & 5\\
Incidence graph of $\mathrm{PG}(3,3)$
                    & 80     & 3        & 14\\
  \hline
  \end{tabular}
\caption{Metric dimension of distance-transitive graphs of valency 13 \label{table:valency13}}
\end{table}

\newpage

\subsection{Strongly regular graphs on up to 45 vertices} \label{subsection:srg}
Recall that a strongly regular graph $\Gamma$ has {\em parameters} $(n,k,a,c)$, where $n$ is the number of vertices, $k$ is the valency, $a$ is the number of common neighbours of a pair of adjacent vertices, and $c$ is the number of common neighbours of a pair of non-adjacent vertices.  In this subsection, we consider the strongly regular graphs on between 35 and 45 vertices, for parameters where a complete classification of graphs are known, and obtain their metric dimension; the results are given in Table~\ref{table:srg}.  (Smaller strongly regular graphs were considered in Tables~\ref{table:smalldrg1} and~\ref{table:smalldrg2} above.)  These graphs were obtained from the online catalogue of Spence~\cite{spence}.  This forms an independent verification of the earlier calculations of Kratica {\em et al.}~\cite{Kratica08}; however, we do not consider strongly regular graphs with parameters $(37,18,8,9)$, as it is unknown if the $6760$ known graphs form the complete set.  (The Paley graphs on 37 and 41 vertices are considered in the next subsection.)

\begin{table}[hbtp]
  \centering
  \begin{tabular}{|l|c|c|c|}
  \hline
  Parameters & No.\ of graphs & Metric dimension & Notes \\
  \hline
  $(35,18,9,9)$ & 3854        & 6 & \\ \hline
  $(36,10,4,2)$ & 1           & 7 & $H(2,6)$ \\ \hline
  $(36,14,4,6)$ & 180         & 6 & \\ \hline
  $(36,14,7,4)$ & 1           & 6 & $J(9,2)$ \\ \hline
  $(36,15,6,6)$ & 32,548      & 6 & \\ \hline
  \multirow{2}{*}{$(40,12,2,4)$}
                & 27          & 7 & \\
                & 1           & 8 & point graph of dual $GQ(3,3)$ \\ \hline
  \multirow{2}{*}{$(45,12,3,3)$}
                & 57          & 7 & \\
                & 21          & 8 & \\ \hline
  $(45,16,8,4)$ & 1           & 7 & $J(10,2)$ \\
  \hline
  \end{tabular}
\caption{Metric dimension of strongly regular graphs on up to 45 vertices \label{table:srg}}
\end{table}

\subsection{Rank-3 strongly regular graphs with up to 100 vertices} \label{subsection:rank3}
A graph which is both strongly regular and distance-transitive is called a {\em rank-$3$ graph}, as its automorphism group has permutation rank $3$ (see \cite{Cameron99}).  Other than a complete multipartite graph, such a graph necessarily has a primitive automorphism group, and as {\sf GAP} contains a library of all primitive groups on up to 2499 points (as obtained by Roney-Dougal~\cite{RoneyD05}), it is straightforward to construct rank-$3$ graphs with a (relatively) small number of vertices.  For this paper, we considered the rank-$3$ graphs with up to $100$ vertices; the metric dimension of each of these graphs is given in Table~\ref{table:rank3}.  (For the Higman--Sims graph on 100 vertices, the exact value was not determined, but an upper bound of $14$ was obtained; the author suspects this is the exact value.)

\begin{table}[hbtp]
  \centering
  \begin{tabular}{|l|c|c|c|}
  \hline
  Graph                     & Parameters      & Metric dimension \\
  \hline
  Paley graph $P_{37}$      & $(37,18,8,9)$   & 5 \\
  Paley graph $P_{41}$      & $(41,20,9,10)$  & 7 \\
  $H(2,7)$                  & $(49,12,5,2)$   & 8 \\
  Paley graph $P_{49}$      & $(49,24,11,12)$ & 7 \\
  Self-complementary graph from $7^2:(3 \times D_{16})$
                            & $(49,24,11,12)$ & 7 \\
  Hoffman--Singleton graph  & $(50,7,0,1)$    & 11 \\
  Paley graph $P_{53}$      & $(53,26,12,13)$ & 7 \\
  $J(11,2)$                 & $(55,18,9,4)$   & 8 \\
  Gewirtz graph             & $(56,10,0,2)$   & 9 \\
  Paley graph $P_{61}$      & $(61,30,14,15)$ & 7 \\
  (from $\mathrm{PSp}(6,2)$)
                            & $(63,30,13,15)$ & 6 \\
  $H(2,8)$                  & $(64,14,6,2)$   & 10 \\
  (from $2^6:(3.S_6)$)      & $(64,18,2,6)$   & 10 \\
  (from $2^6:(S_3 \times \mathrm{GL}(3,2)$)
                            & $(64,21,8,6)$   & 9 \\
  Affine polar graph $\mathrm{VO}^-(6,2)$
                            & $(64,27,10,12)$ & 7 \\
  (from $2^6:S_8$)          & $(64,28,12,12)$ & 7 \\
  $J(12,2)$                 & $(66,20,10,4)$  & 8 \\
  Paley graph $P_{73}$      & $(73,36,17,18)$ & 7 \\
  $M_{22}$ graph            & $(77,16,0,4)$   & 11 \\
  $J(13,2)$                 & $(78,22,11,4)$  & 9 \\
  $H(2,9)$                  & $(81,16,7,2)$   & 11 \\
  Brouwer--Haemers (from $3^4:2.\mathrm{P\Gamma L}(2,9)$) \cite{BrouwerHaemers92}
                            & $(81,20,1,6)$   & 11 \\
  (from $3^4:\mathrm{GL}(2,3):S_4$)
                            & $(81,32,13,12)$ & 8 \\
  Paley graph $P_{81}$      & $(81,40,19,20)$ & 7 \\
  Self-complementary graph from $3^5:(4.S_5)$
                            & $(81,40,19,20)$ & 8 \\
  (from $PSp(4,4).2$)       & $(85,20,3,5)$   & 12 \\
  Paley graph $P_{89}$      & $(89,44,21,22)$ & 8 \\
  $J(14,2)$                 & $(91,24,12,4)$  & 10 \\
  Paley graph $P_{97}$      & $(97,48,23,24)$ & 8 \\
  $H(2,10)$                 & $(100,18,8,2)$  & 12 \\
  Higman--Sims graph        & $(100,22,0,6)$  & $\leq 14$ \\
  Hall--Janko graph         & $(100,36,14,12)$ & 9 \\
  \hline
  \end{tabular}
\caption{Metric dimension of rank-$3$ strongly regular graphs on up to $100$ vertices \label{table:rank3}}
\end{table}

\newpage

\subsection{Hadamard graphs} \label{subsection:Hadamard}

A {\em Hadamard matrix} of order $k$ is a $k\times k$ real matrix $H$ with entries $\pm 1$ and the property that $H H^t = kI$.  Such a matrix must have order $1$, $2$ or a multiple of $4$ (it is conjectured that all multiples of $4$ are admissible).  From a Hadamard matrix $H$ of order $k\geq 4$, the associated {\em Hadamard graph} $\Gamma(H)$ has $4k$ vertices $\{r_i^+, r_i^-, c_i^+, c_i^- \, : \, 1\leq i \leq k\}$, with $r_i^+ \sim c_i^+$ and $r_i^- \sim c_i^-$ if $H_{ij}=1$, and $r_i^+ \sim c_i^-$ and $r_i^- \sim c_i^+$ if $H_{ij}=-1$.  This graph is bipartite and distance-regular, with diameter~$4$ and valency $k$; for further details, see \cite[Section 1.8]{BCN}.

Using Sloane's library of Hadamard matrices, one can easily construct the corresponding Hadamard graphs in $\textsf{GRAPE}$.  The metric dimension for the Hadamard graphs arising from Hadamard matrices of orders from $4$ to $20$ is given in Table~\ref{table:hadamard} below.  In each case considered, Hadamard graphs of the same order had the same metric dimension.

\begin{table}[hbtp]
  \centering
  \begin{tabular}{|c|c|c|c|}
  \hline
  Order ($=$valency) & No. of vertices & No. of examples & Met.\ dim. \\
  \hline
  4                  & 16              & 1 ($\cong Q_4$) & 4 \\
  8                  & 32              & 1               & 7 \\
  12                 & 48              & 1               & 8 \\
  16                 & 64              & 5               & 10 \\
  20                 & 80              & 3               & 10 \\
  \hline
  \end{tabular}
\caption{Metric dimension of Hadamard graphs on up to $80$ vertices \label{table:hadamard}}
\end{table}

\renewcommand{\arraystretch}{1.0}

%\newpage
\subsection*{Acknowledgements}
Most of this work was carried out while the author was a postdoctoral fellow at Ryerson University.  The author would like to thank Leonard Soicher and Alexander Hulpke for their assistance with \textsf{GAP}, and the Faculty of Science, Ryerson University for financial support.

\end{document}